\documentclass[12pt]{article}

\usepackage{amsfonts}

\newtheorem{thm}{Theorem}[section]

\newtheorem{cor}[thm]{Corollary}

\newtheorem{conjecture}{Conjecture}

\newcommand{\be}{\begin{equation}}
\newcommand{\ee}{\end{equation}}
\newcommand{\openbox}{\leavevmode
  \hbox to8pt{\hfil\vrule\vbox to6pt{\hrule width6pt\vfil\hrule}\vrule}}

\newcommand{\qed}{\hbox to5pt{ } \hfill \openbox\bigskip\medskip}

\newcommand{\rk}{\mbox{\rm rank}}

\newcommand{\ve}[1]{\mathbf{#1}}

\newcommand{\cT}{\mbox{$\cal T$}}

\newcommand{\cF}{\mbox{$\cal F$}}
\newcommand{\cG}{\mbox{$\cal G$}}

\newcommand{\cM}{\mbox{$\cal M$}}

\newcommand{\R}{\mathbb R}

\newcommand{\F}{\mathbb F}

\title{A new upper bound for the size of a sunflower-free family}
\author{G\'abor Heged\"{u}s
\\{\normalsize }
}

\begin{document}

\maketitle
\begin{abstract}
We combine here Tao's slice-rank bounding method and Gr\"obner basis techniques and apply here to the Erd\H{o}s-Rado Sunflower Conjecture.

Let $\frac{3k}{2}\leq n\leq 3k$ be integers. We prove that if $\mbox{$\cal F$}$ be a $k$-uniform  family of subsets of  $[n]$ without a sunflower with 3 petals, then
$$
|\mbox{$\cal F$}|\leq  3{n \choose n/3}.
$$ 

We give also some new upper bounds for the size of a sunflower-free family in  $2^{[n]}$.
\end{abstract}


\medskip
\section{Introduction}

First we introduce some notation. 

Let  $[n]$ stand for the set $\{1,2,
\ldots, n\}$. We denote the family of all subsets of $[n]$  by $2^{[n]}$. 

Let $X$ be a fixed subset of $[n]$. For an integer $0\leq k\leq n$ we denote by
${X \choose k}$ the family of all  $k$ element subsets of $X$. This is the {\em complete} $k$-uniform family.

We say that a family $\cF$ is {\em $k$-uniform}, if $|F|=k$ for each $F\in \cF$.

A family $\cF=\{F_1,\ldots ,F_m\}$ of subsets of $[n]$ is a {\em sunflower} (or a {\em $\Delta$-system}) with $t$ petals if
$$
F_i\cap F_j=\bigcap\limits_{s=1}^t F_s
$$
for each $1\leq i,j\leq t$.

Here the intersection of the members of a sunflower form its {\em kernel}.

Erd\H{o}s and  Rado conjectured  the following famous statement in \cite{ER}.
\begin{conjecture} \label{Econj}
For each $t>2$, there exists a constant $C(t)$ such that  if  $\cF$ is a $k$-uniform set system with more than 
$C(t)^k$
members,  then $\cF$ contains a sunflower with $t$ petals.
\end{conjecture}
Erd\H{o}s  offered 1000 dollars for the proof or disproof of this conjecture for $t=3$ (see \cite{E}).

Erd\H{o}s and  Rado gave also an upper bound for the size of a $k$-uniform family without a sunflower with $t$ petals in \cite{ER}.
\begin{thm} \label{Sthm} (Sunflower theorem)
If $\cF$ is a $k$-uniform set system with more than 
$$
k!(t-1)^k\big( 1-\sum_{s=1}^{k-1} \frac{s}{(s+1)! (t-1)^s}\big)
$$ 
members, then $\cF$ contains a sunflower with $t$ petals.
\end{thm}

Define $F(n,t)$ to be the largest integer so that there exists a family $\cF$ of subsets of $[n]$ which does not contain a sunflower with $t$ petals and $|\cF|=F(n,t)$.

Define $\beta_t$ as
$$
\beta_t:= \lim_{n \to \infty} F(n,t)^{1/n}.
$$
Naslund and Sawin gave the following upper bound for the size of a sunflower-free family in \cite{NS}. Their proof based on Tao's slice--rank bounding method (see the blog \cite{T}). 

\begin{thm} \label{NS}
Let $\cF$ be a family of subsets of  $[n]$ without a sunflower with 3 petals. Then
$$
|\cF|\leq 3n\Big( \sum_{i=0}^{n/3} {n\choose i}\Big).
$$ 
\end{thm}

Naslund and Sawin proved also the following upper bound for $\beta_3$ in \cite{NS}.

\begin{cor} \label{NS2}
$$
\beta_3\leq \frac{3}{2^{2/3}}=1.88988...
$$
\end{cor}

Our main result is the following new upper bound for the size of a sunflower-free family. In the proof we mix Tao's slice--rank bounding method with Gr\"obner basis techniques. Our proof is a simple modification of the proof of Theorem 1 in \cite{NS}.

\begin{thm} \label{main-uni}
Let $\frac{3k}{2}\leq n\leq 3k$ be integers. Let $\cF$ be a $k$-uniform  family of subsets of  $[n]$ without a sunflower with 3 petals. Then
$$
|\cF|\leq  3{n \choose n/3}.
$$ 
\end{thm}

Theorem \ref{main-uni} implies easily the following Corollary.

\begin{cor} \label{main}
 Let $\cF$ be a  sunflower--free family of subsets of  $[n]$. Then
$$
|\cF|\leq 3\lceil\frac{n}{3}\rceil {n \choose n/3}+2 \sum_{i=0}^{\lceil n/3 \rceil} {n\choose i}.
$$ 
\end{cor}

In Section 2 we collected some useful preliminaries about the slice rang of functions and Gr\"obner bases.
We present our proofs in Section 3.

\section{Preliminaries}

\subsection{Slice rang}

Let $\delta$ denote in this Section the delta function.

We define first the slice rang of functions. This definition appeared first in Tao's blog \cite{T}. 

Let $A$ be a fixed finite set, $m\geq 1$ be a fixed integer and $\F$ be a field.

Recall that a function $F:A^m\to \F$ has {\em slice-rank} one,  if it has the form 
$$
(\ve x_1, \ldots ,\ve x_m) \mapsto f(\ve x_i)g(\ve x_1, \ldots ,\ve x_{i-1}, \ve x_{i+1}, \ldots ,\ve x_m),
$$
for some $i=1, \ldots ,m$ and some functions $f:A\to \F$, $g:A^{m-1}\to \F$.

The slice rank $\rk(F)$ of a function $F:A^m\to \F$ is the least number of rank one functions needed to generate $F$ as a linear combination.

For instance, if $m=2$, then we get back the usual definition of the rank of a function $F:A^2\to \F$.

Tao  proved the following  result about the slice rang of diagonal hyper-matrices in \cite{T} Lemma 1 (see also Lemma 4.7 in \cite{BCCGU}).
\begin{thm}\label{delta}
Let $\F$ be a fixed field, let $\cT\subseteq {\F}^n$ be a finite subset and let $c_{\alpha}\in \F$ denote a coefficient for each  $\alpha \in \cT$ . Consider the function
$$
F(\ve x_1, \ldots ,\ve x_m):= \sum_{\alpha\in \cT} c_{\alpha}\delta_{\alpha}(\ve x_1)\ldots  \delta_{\alpha}(\ve x_m):{\cT}^m\to \F.
$$
Then 
$$
\rk(F)=|\{\alpha \in \cT:~ c_{\alpha}\ne 0\}|.
$$
\end{thm}

\subsection{Gr\"obner  theory}

Let $\F$ be a field.
In the following $\F[x_1, \ldots, x_n]=\F[\ve x]$ denotes  the
ring of polynomials in commuting variables $x_1, \ldots, x_n$ over $\F$.
For a subset $F \subseteq [n]$ we write
$\ve x_F = \prod_{j \in F} x_j$.
In particular, $\ve x_{\emptyset}= 1$.

We denote by $\ve v_F\in \{0,1\}^n$ the
characteristic vector of a set
$F \subseteq [n]$.
For a family of subsets $\cF \subseteq 2^{[n]}$, define
$V(\cF)$ as the subset  $\{\ve v_F : F \in \cF\} \subseteq \{0,1\}^n \subseteq \F^n$.
A polynomial $f\in \F[x_1,\ldots ,x_n]$ can be considered as a function
from $V(\cF)$ to $\F$ in a natural way.

We can describe several interesting properties of finite set systems ${\mathcal
F}\subseteq 2^{[n]}$  as statements
about  {\em polynomial functions on $V(\cF)$}. 
As for polynomial functions on $V(\cF)$,
it is natural to consider the ideal $I(V(\cF))$:
$$ 
I(V(\cF)):=\{f\in \F[\ve x]:~f(\ve v)=0 \mbox{ whenever } \ve v\in V(\cF)\}. 
$$
Clearly the substitution gives  an $\F$ algebra homomorphism from $\F[\ve x]$ to the $\F$ algebra  of $\F$-valued functions on $V(\cF)$. 
It is easy to verify that this homomorphism is
surjective, and the kernel is exactly
$I(V(\cF))$. Hence we can identify the algebra $\F[\ve x]/I(V(\cF))$ and the algebra of
$\F$ valued functions on $V(\cF)$. It follows that 
$$
\dim _\F
\F[\ve x]/I(V(\cF))=|\cF|.
$$


Now we recall some basic facts about to Gr\"obner 
bases and standard monomials. For details we refer to \cite{AL}, \cite{Bu}, \cite{CCS}, \cite{CLS}.

A linear order $\prec$ on the monomials over 
variables $x_1,x_2,\ldots, x_m$ is a {\em term order}, or {\em monomial
order}, if 1 is
the minimal element of $\prec$, and $\ve u \ve w\prec \ve v\ve w$ holds for any monomials
$\ve u,\ve v,\ve w$ with $\ve u\prec \ve v$. Two important term orders are the lexicographic
order $\prec_l$ and the deglex order $\prec _{d}$. We have
$$x_1^{i_1}x_2^{i_2}\cdots x_m^{i_m}\prec_l x_1^{j_1}x_2^{j_2}\cdots
x_m^{j_m}$$
iff $i_k<j_k$ holds for the smallest index $k$ such
that $i_k\not=j_k$. The definition of  the deglex order is similar: we have $\ve u\prec_{d} \ve v$ iff  either
$\deg \ve u <\deg \ve v$, or $\deg \ve  u =\deg \ve v$, and $\ve u\prec_l \ve v$.

The {\em leading monomial} ${\rm lm}(f)$
of a nonzero polynomial $f\in \F[\ve x]$ is the $\prec$-largest
monomial which appears with nonzero coefficient in the canonical form of $f$ as a linear
combination of monomials.

Let $I$ be an ideal of $\F[\ve x]$. We say that a finite subset $\cG\subseteq I$ is a {\it
Gr\"obner basis} of $I$ if for every $f\in I$ there exists a $g\in \cG$ such
that ${\rm lm}(g)$ divides ${\rm lm}(f)$. In other words, the leading
monomials ${\rm lm}(g)$ for $g\in \cG $ generate the  ideal of 
monomials
$\{ {\rm lm}(f):~f\in I\}$. Consequently $\cG$ is
actually a basis of $I$, i.e. $\cG$ generates $I$ as an ideal of $\F[\ve x]$. 
A well--known fact is (cf. \cite[Chapter 1, Corollary
3.12]{CCS} or \cite[Corollary 1.6.5, Theorem 1.9.1]{AL}) that every
nonzero ideal $I$ of $\F[\ve x]$ has a Gr\"obner basis.

A monomial $\ve w\in \F[\ve x]$ is a {\it standard monomial for $I$} if
it is not a leading monomial for any $f\in I$. We denote by  ${\rm sm}(I)$ 
the set of standard monomials of $I$.

Let $\cF\subseteq 2^{[n]}$ be a set family. Then   the characteristic vectors 
in $V(\cF)$ are all 0,1-vectors, consequently the polynomials $x_i^2-x_i$ 
all vanish  on $V(\cF)$. It follows that the standard monomials of the ideal
$I(\cF):=I(V(\cF))$ are square-free monomials.

Now we give a short introduction to the notion of reduction.
Let $\cG$ be a set of polynomials in $\F[x_1,\ldots,x_n]$ and
let $f\in \F[x_1,\ldots,x_n]$ be a fixed polynomial.
We can reduce $f$ by the set $\cG$ with respect to $\prec$.
This gives a new polynomial $h \in \F[x_1,\ldots,x_n]$.

Here {\em reduction} means that we possibly repeatedly replace monomials
in $f$ by smaller ones (with respect to $\prec$) in the following way: if $w$ is a
monomial occurring in $f$ and ${\rm lm}(g)$ divides $w$ for some
$g\in {\cal G}$ (i.e. $w={\rm lm}(g)u$ for some monomial $u$), then
we replace $w$ in $f$ with $u({\rm lm}(g)-g)$. It is easy to verify that the
monomials in $u({\rm lm}(g)-g)$
are $\prec$-smaller than $w$.


It is a key fact that
${\rm sm}(I)$ gives a basis of the $\F$-vector-space $\F[\ve x]/I$ in
the sense that  
every polynomial $g\in \F[\ve x]$ can be 
uniquely expressed as $h+f$ where $f\in I$ and
$h$ is a unique $\F$-linear combination of monomials from ${\rm sm}(I)$. Hence if $g\in \F[\ve x]$ is an arbitrary polynomial and $\cG$ is a Gr\"obner basis of $I$, then we can reduce $g$ with $\cG$ into a linear combination of
standard monomials for $I$. In particular, $f\in I$ if and only if $f$ can be 
$\cG$-reduced to 0.  

Let $0\leq k\leq n/2$ and denote by $\cM_{k,n}$ the set of all
monomials $\ve x_G$ such that $G=\{s_1<s_2<\ldots <s_j\}\subset [n]$
for which $j\leq k$ and $s_i\geq 2i$ holds for every $i$, $1\leq i\leq j$. These monomials $\ve x_G$ are the {\em ballot monomials} of degree at most $k$.
If $n$ is clear from the context, then we write $\cM_k$ instead of the more precise $\cM_{k,n}$. 
 It is known 
that 
$$ 
|\cM_k|={n \choose k}.
$$

In \cite{HR} we described completely the Gr\"obner bases and the standard monomials  of the complete uniform families of all $k$ element subsets of $[n]$.

\begin{thm}\label{sm1}
Let $\prec$ an arbitrary term order such that $x_1\prec \ldots \prec x_n$. Let $0\leq k\leq n$ and  
$j:=min(k,n-k)$. Then 
$$                                
{\rm sm}(V{[n]\choose k})=\cM_{j,n}.
$$
\end{thm}

Let $0\leq k\leq n$ and $\ell>0$ be arbitrary integers. Define the vector system
$$
\cF(n,k,\ell):=\underbrace{V({[n]\choose k})\times \ldots \times V({[n]\choose k})}_{\ell}\subseteq \{0,1\}^{n\ell}.
$$

It is easy to verify the following Corollary.

\begin{cor} \label{sm2}
Let $\prec$ an arbitrary term order such that $x_1\prec \ldots \prec x_n$. Let $0\leq k\leq n$ and $\ell>0$ be arbitrary integers. Let 
$j:=min(k,n-k)$. Then 
$$                                
{\rm sm}(\cF(n,k,\ell))=\{x_{M_1}\cdot \ldots \cdot x_{M_{\ell}}:~ x_{M_1},\ldots ,x_{M_{\ell}}\in\cM_{j,n} \}.
$$
\end{cor}

\section{Proofs}

{\bf Proof of Theorem \ref{main-uni}:} 

Let $\cF$ be a $k$-uniform sunflower--free family of subsets of  $[n]$.

Let $H_1,H_2,H_3 \in  \cF$ be arbitrary subsets. Since $\cF$ is sunflower--free, hence if 
$$
\ve v(H_1)+\ve v(H_2)+\ve v(H_3)\in \{0,1,3\}^n,
$$
then $H_1=H_2=H_3$.

Namely first suppose that $H_1\ne H_2$, $H_1\ne H_3$ and $H_2\ne H_3$. Then the triple $(H_1,H_2,H_3)$ is not a sunflower, hence there exist indices $1\leq i<j\leq 3$ such that $(H_i\cap H_j)\setminus (H_1\cap H_2\cap H_3)\ne \emptyset$. Let  $t\in (H_i\cap H_j)\setminus (H_1\cap H_2\cap H_3)$. Then $\ve v(H_1)_t+\ve v(H_2)_t+\ve v(H_3)_t =2$. 

Suppose that $H_1\ne H_2$ but $H_2=H_3$. Since $|H_1|=|H_2|=k$, hence $H_2\setminus H_1\ne \emptyset$. Let $t\in H_2\setminus H_1$. Then it is easy to see that $\ve v(H_1)_t+\ve v(H_2)_t+\ve v(H_3)_t =2$. 

Consider the polynomial function 
$$
T:(\cF)^3\to \R
$$
given by
$$
T(\ve x,\ve y,\ve z):=\prod_{i=1}^n (2-(x_i+y_i+z_i))
$$
for each $\ve x=(x_1,\ldots ,x_n), \ve y=(y_1,\ldots ,y_n),\ve z=(z_1,\ldots ,z_n)\in \cF\subseteq V{[n]\choose k}$.

Let $\cG$ denote a deglex  Gr\"obner basis of the ideal  $I:=I(\cF(n,k,3))$. Let $H$ denote the reduction of $T$ via $\cG$.

Then 
\be \label{HTeq}
H(\ve x,\ve y,\ve z)=T(\ve x, \ve y,\ve z)
\ee
for each  $\ve x=(x_1,\ldots ,x_n), \ve y=(y_1,\ldots ,y_n),\ve z=(z_1,\ldots ,z_n)\in \cF\subseteq V{[n]\choose k}$, because we reduced $T$ with a Gr\"obner basis of the ideal $I$.

On the other hand 
$$
T(\ve x,\ve y,\ve z)\ne 0 \mbox{ if and only if }\ve x=\ve y=\ve z\in \cF,
$$
hence by equation (\ref{HTeq})
\be \label{delta2}
H(\ve x,\ve y,\ve z)\ne 0 \mbox{ if and only if }\ve x=\ve y=\ve z\in \cF.
\ee

Let $j:=min(k,n-k)$.

Since $\cF\subseteq V{[n]\choose k}$, hence it follows from Corollary \ref{sm2} that we can write $H(\ve x,\ve y,\ve z)$ as a linear combination of standard monomials 
$$
x_Iy_Kz_L,
$$
where $x_I,y_K,z_L\in \cM_{j,n}$ and $deg(x_Iy_Kz_L)\leq n$. Here we used that $\cG$ is a {\em deglex} Gr\"obner basis of the ideal $I$.

It follows from the pigeonhole principle that at least one of $|I|$, $|K|$ and $|L|$ is at most $n/3$. 

First we can consider the contribution of the standard  monomials to the sum for which $|I|\leq \frac{n}{3}$. 

We can regroup this contribution as 
$$
\sum_{M} x_Mg_{M}(\ve y,\ve z),
$$
where $M$ ranges over those subsets $\{i_1, \ldots ,i_t\}$ of $[n]$ with $t\leq n/3$ and $i_s\geq 2s$ for every $1\leq s\leq t$. Here $g_{M}:(\cF)^2\to \R$ are some explicitly computable functions.

The number of such $M$ is at most ${n\choose n/3}$, so this contribution has slice--rank at most 
${n\choose n/3}$. 

The remaining contributions arising from the cases $|K|\leq \frac{n}{3}$ and  $|L|\leq \frac{n}{3}$. 

But $H$ and $T$ are the same functions on $\cF(n,k,3)$, hence we get that
$$
\rk(H)=\rk(T)\leq 3 {n\choose n/3}.
$$

But it follows from Theorem \ref{delta} and  (\ref{delta2}) that
$$
\rk(H)=|\cF|. 
$$

Finally we get that
$$
|\cF|\leq 3 {n\choose n/3}.
$$
\qed

{\bf Proof of Corollary  \ref{main}:} 
Let $\cF\subseteq \{0,1\}^n$ be a fixed sunflower-free subset.  Define
the families
$$
\cF(s):=\cF \cap {[n]\choose s}
$$
for each $0\leq s\leq n$.

We have two separate cases.
                          
1. Suppose that either $s\geq \frac{2n}{3}$ or $s\leq \frac{n}{3}$. Then clearly
$$
|\cF(s)|\leq {n\choose s}.
$$

2. Suppose that $\frac{n}{3}\leq s\leq \frac{2n}{3}$.

Then we can apply Theorem \ref{main-uni} for the family $\cF(s)$ and we get
$$
|\cF(s)|\leq 3{n \choose n/3}.
$$

Finally 
$$
|\cF|= \sum_{s=0}^n |\cF(s)|\leq \lceil\frac{n}{3}\rceil\Big( 3{n \choose n/3}\Big)+2 \sum_{i=0}^{\lceil n/3\rceil} {n\choose i}.
$$
\qed

\section{Concluding remarks}

It is easy to verify that Conjecture \ref{Econj} follows immediately from the following conjecture.

\begin{conjecture} \label{Hconj}
For each $t>2$  there exists a constant $C(t)$ such that  if  $\cF$ is an arbitrary  $k$-uniform set system, which does not contain any sunflower with $t$ petals, then  $|\cup_{F\in \cF} F|\leq C(t)k$. 
\end{conjecture}


The following Corollary is clear.

\begin{cor} \label{HG}           
Suppose that Conjecture \ref{Hconj} is true for $t=3$. Let $\cF$ be an arbitrary  $k$-uniform set system, which does not contain any sunflower with $3$ petals. Then 
$$
|\cF| \leq 3{C(t)k\choose C(t)k/3}.
$$
\end{cor}

{\bf Acknowledgements.}           
I am very indebted to  Lajos R\'onyai for his useful remarks.

\end{document}